\input amstex

\def\BC{{\Bbb C}}
\def\BN{{\Bbb N}}
\def\BP{{\Bbb P}}
\def\BQ{{\Bbb Q}}
\def\BZ{{\Bbb Z}}

\def\CD{{\Cal D}}
\def\CL{{\Cal L}}
\def\CM{{\Cal M}}
\def\CN{{\Cal N}}
\def\CO{{\Cal O}}

\def\CT{{\Cal T}}
\def\ft{{\frak t}}
\def\tphi{\tilde{\phi}}


\def\btu{\bigtriangleup}
\def\dpar{\partial}
\def\hra{\hookrightarrow}
\def\iso{\buildrel\sim\over\longrightarrow} 

\def\lra{\longrightarrow}

\parskip=6pt

\documentstyle{amsppt}
\document
\magnification=1200

\centerline{\bf ON HYPERGEOMETRIC FUNCTIONS} 
\centerline{\bf CONNECTED 
WITH QUANTUM COHOMOLOGY OF FLAG SPACES}
\bigskip
\centerline{Vadim Schechtman}
\smallskip
\centerline{Max-Planck-Institut f\"ur Mathematik,}
\centerline{Gottfried-Claren-Stra\ss e 26, 53225 Bonn}
\bigskip

\bigskip
\centerline{\bf Introduction}
\bigskip

In the Givental's work on the Gromov-Witten invariants for projective 
complete intersections, [G1], the principal role is played by certain 
formal power series connected with the quantum cohomology of a 
manifold. One has a manifold $X$  
with a natural torus action, with the finite number of fixed points, 
and one has a power series $Z_w$ associated with each fixed point $x_w$. 
The coefficients of these series are certain integrals over the spaces 
of stable maps of genus $0$ curves with two marked points to $X$. 
These series form a fundamental 
system of solutions of a certain lisse $\CD$-module on a power of the 
punctured disk. The small quantum cohomology of $X$ coincides with 
the algebra of functions on its characteristic variety. 
The series $Z_w$ are uniquely determined by certain {\bf recursion 
relations} relating $Z_w$ with all $Z_{w'}$ if $x_w$ is connected with 
$x_{w'}$ by a fixed line.

The present work consists of three parts. The First part contains nothing 
new. Here we present the Givental's computations from [G1], 
in the simplest case of a projective space, in more detail than in [G1]. 
In the Second part, we write down the above recursion relations for the 
flag spaces $X=G/B$ ($G$ being a simple algebraic group), see {\bf Theorem II.3.8}, 
which is the main result of this paper. Here $X$ is equipped with the natural 
action of the maximal torus of $G$. 

In the work [G2], Givental gave another beautiful set of relations which 
also determines the above mentioned series completely. Namely, these 
are Toda lattice differential equations (more precisely, 
we need the equivariant version of the results of [G2]). This set of relations 
has completely different nature, and it is highly non-trivial fact 
that both sets of relations determine the same series. In the Third part, 
we check this by a direct computation for $G=SL(3)$. It turns out  
that in this case the series $Z_w$ admit a nice explicit expression, 
see III.1.2, III.2.2. ({\it A posteriori} it is not surprising, since 
in this case $X$ admits the Pl\"ucker embedding into $\BP^1\times\BP^1$, and one can use another computation by Givental, dealing with the toric complete 
intersections). 

This paper arose from the author's attempt to make an exercise on the 
beautiful fixed point technique of Givental and Kontsevich. I am grateful 
to Yu. I. Manin for many stimulating discussions, to A. Goncharov for the 
stimulating interest on the first stage of this work, and to V. Batyrev 
who communicated to me formula III.1.4. This work was done in the 
highly stimulating atmosphere of Max-Planck-Institut f\"ur Mathematik, 
and I am very grateful to MPI for the hospitality.


\newpage
\bigskip
\centerline{\bf Part I. PROJECTIVE SPACES}
\bigskip
\centerline{\bf \S 1. Equivariant cohomology of $\BP^n$}
\medskip

{\bf 1.1.} Let $X$ denote the $n$-dimensional projective space $\BP^n=
\{(z_0:\ldots:z_n)\}$, the space of lines $L\subset V=\BC^{n+1}$. 
The torus $T=(\BC^*)^{n+1}$ acts on $X$ by the rule
$(\alpha_0,\ldots,\alpha_n)\cdot(z_0:\ldots:z_n)=
(\alpha_0z_0:\ldots:\alpha_nz_n)$; this action has $n+1$ fixed 
points
$x_i=(0:\ldots 0:1:0\ldots 0)$ (one on $i$-th place), $i=0,\ldots,n$. 

Let $\CL$ denote the line bundle over $X$ whose fiber over $L\subset V$ 
is $L$; $\CL$ has an obvious $T$-equivariant structure. 

Let $\ft$ denote the Lie algebra of $T$, $\ft^*_{\BZ}:=Hom (T,\BC^*)
\subset \ft^*$ the lattice of characters. For $\lambda\in\ft^*_{\BZ}$, 
let $L_{\lambda}$ be the $T$-equivariant line bundle over the point, 
equal to $\BC$, with $T$ acting by means of the character $\lambda$. 
Assigning to $\lambda$ the first Chern class $c_1(L_{\lambda})
\in H^2_T(pt)$, we identify the graded ring $A:=H^*_T(pt)$ with $\BC[t_{\BZ}^*]=
\BC[\lambda_0,\ldots,\lambda_n]$, $\lambda_i$ being the projection 
on $i$-th factor. 

The graded $A$-algebra $R:=H^*_T(X)$ is identified with $\BC[p,\lambda_0,\ldots,
\lambda_n]/(\prod_i(p-\lambda_i))$ where $p:=c_1(\CL)\in H^2_T(X)$. It is computed using the Bott's fixed point 
theorem, [AB]. 

Since our cohomologies will be even anyway, it is 
convenient to use the grading of the rings $A$, $R$, etc. 
by assigning to $p$ and $\lambda_i$ the degree $1$.  

{\bf 1.2.} The Euler classes of the tangent spaces 
at the fixed points $x_i$ are equal to 
$$
e_i:=e(T_{X;x_i})=\prod_{b\neq i}\ (\lambda_i-\lambda_b)
\in H^{2n}_T(pt)=A^{n}
\eqno{(1.1)}
$$
Let us consider the restriction map 
$i^*_b:\ H^*_T(X)\lra H^*_T(x_b)$. We have 
$$
i^*_b(p)=i^*_b(c_1(\CL))=c_1(\CL_{x_b})
$$
The fiber $\CL_{x_b}$ is the line $\BC\cdot (0,\ldots,1,
\ldots 0)\subset\BC^{n+1}$ ($1$ on $b$-th place); the Lie 
algebra $\ft$ acts on this line by means of the character 
$\lambda_b$. Therefore, 
$$
i^*_b(p)=\lambda_b, 
$$
whence
$$
i^*_b(f(p))=f(\lambda_b)
\eqno{(1.2)}
$$
Let $A'$ be the ring obtained from $A$ by inverting all 
elements $e_i$, i.e. by inverting all the differences 
$\lambda_a-\lambda_b\ (a\neq b)$. Bott's theorem says 
that the restriction map  
$$
i^*:\ H^*_T(X)\lra H^*_T(X^T)=\oplus_{b=1}^n\ A\cdot 1_{x_b}
$$
becomes an isomorphism after the base change $A\lra A'$. 
We denote $R'=R\otimes_AA'$. 

Let us introduce the elements 
$$
\phi_i(p):=\prod_{b\neq i}\ (p-\lambda_b)\in R^{n}
\eqno{(1.3)}
$$
Obviously, 
$$
\phi_i(\lambda_j)=e_i\delta_{ij}
\eqno{(1.4)}
$$
It follows from the Bott's theorem that the set 
$$
\{\phi_i(p)/e_i^{1/2};\ i=0,\ldots,b\}
\eqno{(1.5)}
$$
is the basis of orthonormal idempotents of the algebra 
$R'$ (to be precise, we should adjoin the square roots 
of $e_i$ to $R'$).  

One can express this in a slightly different way. We have 
the integration map
$$
\int_X:\ R\lra A
$$
of degree $-2n$, given by 
$$
\int_X\ f(p)=\sum_i\ res_{p=\lambda_i}\left(\frac{f(p)}
{\prod_b(p-\lambda_b)}dp\right)=\sum_i\ \frac{f(\lambda_i)}
{e_i}
\eqno{(1.6)}
$$
We have the Poincar\'e pairing  
$\langle\cdot,\cdot\rangle:\ R\times R\lra A$, 
$$
\langle f,g\rangle=\int_X\ fg
$$
We have
$$
\langle f,\phi_i\rangle=f(\lambda_i),
\eqno{(1.7)}
$$
hence
$$
\langle\phi_i,\phi_j\rangle=e_i\delta_{ij}
\eqno{(1.8)}
$$
For each $f\in R$, 
$$
f=\sum_i\ \frac{f(\lambda_i)}{e_i}\phi_i
\eqno{(1.9)}
$$
Note that the rhs {\it a proiori} lies in $R'$ but in 
fact it belongs to $R\subset R'$, since the lhs does.

\bigskip
\centerline{\bf \S 2. Partition function}
\medskip

{\bf 2.1.} Let $X_d\ (d\geq 0)$ denote the stack of stable maps $\{ f:(C;y_1,y_2)
\lra X\}$ of genus $0$ curves with two marked points, such that 
$f_*([C])=d\cdot\beta$ where $\beta\in H_2(X)$ is the generator dual to $p$. It is a Deligne-Mumford stack 
for $d\geq 1$.   

Let $\CL_1$ be the line bundle over $X_d$ whose fiber at a point 
$(f,\ldots)$ is the tangent space $\CT_{C;y_1}$; denote 
$c(d)=c_1(\CL_1)\in H^2_T(X_d)$.   

We want to calculate the following formal power series
$$
Z(q,p)=Z(q,p,\lambda,h)=1+\sum_{d\geq 1}\ e_{1*}\left(\frac{1}{h+c(d)}
\right)q^d
\eqno{(2.1)}
$$
Here 
$$
e_1:\ X_d\lra X
$$
is the evaluation map sending $(f,\ldots)$ to $f(y_1)$. This is the same 
as to compute the series
$$
Z_i(q)=\langle Z(q),\phi_i\rangle=\sum_d\ \int_X\phi_i e_{1*}
\left(\frac{1}{h+c(d)}\right) q^d=
\sum_d\ \left(\int_{X_d}\ \frac{e_1^*\phi_i}{h+c(d)}\right) q^d
\eqno{(2.2)}
$$
$(i=0,\ldots,n)$. First, let us formulate the answer.  

{\bf 2.2.} Define the series
$$
S(q,p)=\sum_{d\geq 0}\ \frac{1}{\prod_{b=0}^n\prod_{m=1}^d
(p-\lambda_b+mh)}\cdot q^d
\eqno{(2.3)}
$$
and
$$
S_i(q):=S(q,\lambda_i)=\sum_{d\geq 0}
\frac{1}{d!\prod_{b\neq i}
\prod_{m=1}^d(\lambda_i-\lambda_b+mh)}\cdot\frac{q^d}{h^d}
\eqno{(2.4)}
$$
{\bf 2.3. Theorem.} {\it (a) $Z(q,p)=S(q,p)$. 

(b) For all $i=0,\ldots,n$, $Z_i(q)=S_i(q)$.} 

Of course, (a) and (b) are equivalent, due to the remarks 
of the previous section. The theorem will be proven 
in Section 5, after preliminaries in Sections 3, 4.

{\bf 2.4. Dimension count.} The dimension of $X_d$ is equal 
to
$$
\dim (X_d)=nd+n+d-1
\eqno{(2.5)}
$$
Indeed, one sees easily that the dimension of the space 
of maps $\BP^1\lra\BP^n$ of degree $d$ is equal to 
$(d+1)(n+1)-1$. To get the dimension of $X_d$, we have 
to subtract from this number $3$ (reparametrizations 
of $\BP^1$) and add $2$ (marked points). 

The theorem says that for each $d\geq 0$, 
$$
e_{1*}\left(\frac{1}{h+c(d)}\right)=
\frac{1}{\prod_{b=0}^n\prod_{m=1}^d(p-\lambda_b+mh)}
\eqno{(2.6)}
$$
Let us assign to the variable $h$ the degree $1$.    
The map $e_{1*}$ decreases the degree 
by $\dim(X_d)-\dim(X)=nd+d-1$. Therefore, the 
degree of the lhs of (2.6) is equal to 
$$
-1-(nd+d-1)=-(n+1)d,
$$
which equals the degree of the rhs. Note that this is 
true for $d=0$ as well (the map $e_{1*}$ has a {\it 
positive} degree $1$ in this case!).

\bigskip

\bigskip
\centerline{\bf \S 3. Recursion relation}
\medskip

In this section we will study the series $S(q,p)$ and $S_i(q)$. 

{\bf 3.1.} As a first remark, note that the {\it normalized} series 
$$
S_i^{norm}:=\frac{\langle S,\phi_i\rangle}
{\langle\phi_i,\phi_i\rangle}
\eqno{(3.1)}
$$
may be written as 
$$
S_i^{norm}(q)=\sum_{d\geq 0}\frac{1}{\prod_{b\neq i}\prod_{m=0}^d
(\lambda_i-\lambda_b+mh)}\cdot\frac{q^d}{d!h^d}=
$$
$$
=\sum_{d\geq 0}\frac{1}{\prod_{b=0}^n\prod_{m=0;(b,m)\neq (i,0)}^d
(\lambda_i-\lambda_b+mh)}\cdot q^d
\eqno{(3.2)}
$$
{\bf 3.2.} Define the series
$$
s_i(q):=S_i(qh)=\sum_{d\geq 0} b_i(d)q^d
\eqno{(3.3)}
$$
where
$$
b_i(d)=b_i(d,\lambda,h)=\frac{1}{d!\prod_{j\neq i}\prod_{m=1}^d
(\lambda_i-\lambda_j+mh)}
\eqno{(3.4)}
$$
The next theorem is the main result about our series. 

{\bf 3.3. Theorem.} {\it For each $i=0,\ldots,n$ have 
$$
s_i(q,h)=1+\sum_{k>0}\sum_{j\neq i}\frac{c_i^j(k)\cdot q^k}
{\lambda_i-\lambda_j+kh}\cdot s_j(q;\frac{\lambda_j-\lambda_i}{k})
\eqno{(3.5)}
$$
where
$$
c_i^j(k)=\frac{1}{k!\prod_{b\neq i}\prod_{m=1;\ (b,m)\neq (j,k)}^k
\left(\lambda_i-\lambda_b+\frac{m(\lambda_j-\lambda_i)}{k}\right)}
\eqno{(3.6)}
$$
The relations (3.5), (3.6) uniquely determine the series 
$s_i$.}

The theorem is a variant of a simple fractions 
decomposition. We use the following elementary fact. 

{\bf 3.4. Lemma.} {\it Let $f(h)$ be a non-constant 
polynomial with distinct roots 
\newline $\alpha_1, \ldots, \alpha_N$. 
Then 
$$
\frac{1}{f(h)}=\sum_{k=1}^N\ 
\frac{1}{(h-\alpha_k)f'(\alpha_k)}
\eqno{(3.7)}
$$
We have
$$
f'(\alpha_k)=\left(\frac{f(h)}{h-\alpha_k}\right)
\biggl|_{h=\alpha_k}
\eqno{(3.8)}
$$}

Indeed, the difference of the rhs and lhs of (3.7) does not 
have singularities on $h\in\BP^1$, hence it is a constant; 
but the value of both sides at $\infty$ is equal to $0$, 
hence they are equal. The formula (3.8) is evident. $\btu$

Now let us apply this to the coefficients $b_i(d)$, (3.4);  
we get
$$
b_i(d)=\frac{1}{d!}\sum_{j\neq i}
$$
$$
\sum_{k=1}^d
\frac{1}{\lambda_i-\lambda_j+kh}\cdot
\frac{1}{\prod_{b\neq i}\prod_{m=1;\ (b,m)\neq (j,k)}^d
(\lambda_i-\lambda_b+m\cdot\frac{\lambda_j-\lambda_i}{k})}
\eqno{(3.9)}
$$
Let us split the product in the denominator in the two 
parts: 
$$
\prod_{m=1}^d=\prod_{m=1}^k\cdot\prod_{m=k+1}^d
$$
The second product is equal to 
$$
\prod_{m=k+1}^d:=\prod_{b\neq i}\prod_{m=k+1}^d
(\lambda_i-\lambda_b+m\frac{\lambda_j-\lambda_i}{k})=
$$
$$
=\prod_{b\neq i}\prod_{m'=1}^{d-k}(\lambda_i-\lambda_b+
\lambda_j-\lambda_i+m'\cdot\frac{\lambda_j-\lambda_i}{k})=
\prod_{b\neq i}\prod_{m'=1}^{d-k}
(\lambda_j-\lambda_b+m'\cdot\frac{\lambda_j-\lambda_i}{k})
$$
Hence
$$
\frac{1}{\prod_{m=k+1}^d}=\frac{d!}{k!}b_j(d-k,
\frac{\lambda_j-\lambda_i}{k})
\eqno{(3.10)}
$$
Therefore, 
$$
b_i(d,h)=\sum_{j\neq i}\sum_{k=1}^d\frac{1}{\lambda_i-
\lambda_j+kh}\cdot
\frac{b_j(d-k,\frac{\lambda_j-\lambda_i}{k})}{k!\prod_{b\neq i}
\prod_{m=1;\ (b,m)\neq (j,k)}^k(\lambda_i-\lambda_b+m\cdot\frac{\lambda_j
-\lambda_i}{k})}=
$$
$$
=\sum_{j\neq i}\sum_{k=1}^d\frac{1}{\lambda_i-\lambda_j+kh}\cdot c_i^j(k)\cdot
b_j(d-k,\frac{\lambda_j-\lambda_i}{k})
\eqno{(3.11)}
$$
Obviously, (3.11) is equivalent to (3.5). 
This proves (3.5). 

The uniqueness is obvious since $s_i(0)=1$, and 
the recursion relations determine  
$s_i(q)$ modulo $q^{k+1}$ once we know $s_j(q)$ modulo 
$q^k$ for all $j\neq i$.       

The theorem is proved. $\btu$

{\bf 3.5.} In order to get a better feeling what is going 
on, let us consider some examples. First, {\bf the case 
$n=0$} is almost trivial, but gives a nice answer. 

We have in this case $A=\BC[\lambda]$ 
($\lambda:=\lambda_0$); $R=\BC[p,\lambda]/(p-\lambda)=A$; 
$$
S(q)=S_0(q)=e^{q/h}
\eqno{(3.12)}
$$
There are no recursion relations. 

{\bf 3.6. The case $n=1$.} We have 
$$
R=\BC[p,\lambda_0,\lambda_1]/((p-\lambda_0)(p-\lambda_1));
$$
$$ 
\phi_0=p-\lambda_1,\ \phi_1=p-\lambda_0;\ 
e_0=\lambda_0-\lambda_1,\ e_1=\lambda_1-\lambda_0. 
$$
It is convenient to introduce the "root" 
$\alpha:=\lambda_0-\lambda_1$. We have 
$$
s_0(q)=s_0(q;\alpha)=
\sum_{d\geq 0}\frac{q^d}{d!\prod_{m=1}^d(\alpha+mh)};\ 
\eqno{(3.13)}
$$
$$
s_1(q)=s_1(q;\alpha)=
\sum_{d\geq 0}\frac{q^d}{d!\prod_{m=1}^d(-\alpha+mh)}=
s_0(q;-\alpha)
\eqno{(3.13)'}
$$
The recursion relations look as follows
$$
s_0(q;\alpha,h)=1+\sum_{k>0}\frac{c(k;\alpha)q^k}
{\alpha+kh}\cdot s_1(q;\alpha;-\frac{\alpha}{k})
\eqno{(3.14)}
$$
and
$$
s_1(q;\alpha,h)=1+\sum_{k>0}\frac{c(k;-\alpha)q^k}
{-\alpha+kh}\cdot s_0(q;\alpha;\frac{\alpha}{k})
\eqno{(3.14)'}
$$
where
$$
c(k;\alpha)=\frac{k^k}{(k!)^2}\cdot\frac{1}{\alpha^{k-1}}
\eqno{(3.15)}
$$
The first few values of $c(k;\alpha)$:
$$
c(1;\alpha)=1;\ c(2;\alpha)=\frac{1}{\alpha};\ 
c(3;\alpha)=\frac{3}{4\alpha^2}
\eqno{(3.16)}
$$
The relation $(3.14)'$ is obtained from $(3.14)$ 
by switching $\alpha$ to $-\alpha$. 

Now let us make a little computation: start building up 
the series $s_0, s_1$ using $(3.14)$, $(3.14)'$. We have
$$
s_0(h)=1+\frac{q}{\alpha+h}s_1(-\alpha)+
\frac{q^2}{\alpha+2h}\cdot\frac{1}{\alpha}s_1(-\frac
{\alpha}{2})+\ldots
$$
$$
s_1(h)=1+\frac{q}{-\alpha+h}s_0(\alpha)+\ldots
$$
Thus, 
$$
s_0=1+\frac{q}{\alpha+h}+\ldots;\ 
s_1=1+\frac{q}{-\alpha+h}+\ldots
$$
hence 
$$
s_0=1+\frac{q}{\alpha+h}+\left(\frac{1}{\alpha+2h}
\cdot\frac{1}{\alpha}+\frac{1}{\alpha+h}\cdot\frac{1}
{-2\alpha}\right)q^2=
$$
$$
=1+\frac{q}{\alpha+h}+\frac{q^2}{2(\alpha+h)(\alpha+2h)}+
\ldots
\eqno{(3.17)}
$$
which is the correct answer, up to this order. 

{\bf 3.7.} As the last example, assume that $n$ 
is arbitrary and let us check (3.3), (3.4) 
up to the {\it first} order({\it sic!}) using (3.5), 
(3.6). 

We have
$$
c_i^j(1)=\frac{1}{\prod_{b\neq i,j}(\lambda_j-\lambda_b)}
\eqno{(3.18)}
$$
Therefore, 
$$
s_i(h)=1+\sum_{j\neq i}\frac{q}{\lambda_i-\lambda_j+h}
\cdot\frac{1}{\prod_{b\neq i,j}(\lambda_j-\lambda_b)}+
\ldots=
$$
$$
=1+\frac{q}{\prod_{j\neq i}(\lambda_i-\lambda_j+h)}+\ldots
\eqno{(3.19)}
$$
where we have used the formula
$$
\frac{1}{\prod_{j\neq i}(\lambda_i-\lambda_j+h)}=
\sum_{j\neq i}\frac{1}{\prod_{b\neq i,j}
(\lambda_j-\lambda_b)}\cdot\frac{1}{\lambda_i-\lambda_j+h}
\eqno{(3.20)}
$$

\bigskip
\centerline{\bf \S 4. First reduction.}
\medskip

{\bf 4.1.} We start proving Theorem 2.3.  

Let us define the series $z_i$ by
$$
z_i(q)=Z_i(qh)
\eqno{(4.1)}
$$
(cf. (3.3)). According to Theorem 3.3, in order to prove 
Theorem 2.3, it suffices to show that $z_i(q)$ satisfy the relations 
(3.5). 

{\bf 4.2.} Let us define the coefficients $B_i$ by
$$
S_i(q)=\sum_d\ B_i(d)q^d
\eqno{(4.2)}
$$
(cf. (2.4)). Thus, we have
$$
B_i(d)=\frac{b_i(d)}{h^d}=\frac{1}{\prod_{b=0}^n\prod_{m=1}^d
(\lambda_i-\lambda_b+mh)}
\eqno{(4.3)}
$$
As we have already noted, (3.5) is equivalent to the identities
$$
b_i(d,h)=\sum_{j\neq i}\sum_{k=1}^d\frac{c_i^j(k)}
{\lambda_i-\lambda_j+kh}\cdot b_j(d-k,\frac{\lambda_j-\lambda_i}{k})
\eqno{(4.4)}
$$
or
$$
b_i(d,h)=\sum_{j\neq i}\sum_{k=1}^d\frac{c_i^j(k)}
{\lambda_i-\lambda_j+kh}\cdot\left(\frac{\lambda_j-\lambda_i}{k}\right)
^{d-k}\cdot 
B_j(d-k,\frac{\lambda_j-\lambda_i}{k})
\eqno{(4.5)}
$$
If we assign to $\lambda_i$ and $h$ the degree $1$ then we have 
$$
deg\ b_i(d)=-dn;\ deg\ B_i(d)=-dn-n;\ deg\ c_i^j(k)=-kn+1
\eqno{(4.6)}
$$
and the identities (4.4), (4.5) are homogeneous, cf. 2.4. 

{\bf 4.3.} Let us denote
$$
U_i(d):=\int_{X_d}\frac{e_1^*\phi_i}{h+c(d)}
\eqno{(4.7)}
$$
Thus, 
$$
Z_i(q)=\sum_{d\geq 0}\ U_i(d)q^d
\eqno{(4.8)}
$$
We have 
$$
U_i(d)=\frac{1}{h}\sum_{a\geq 0}\frac{(-1)^a}{h^a}\int_{X_d}
e_1^*\phi_i\cdot c(d)^a
$$
We have $U_i(0)=1$. Assume now that $d\geq 1$. The degree of 
the integral 
\newline $\int_{X_d}e_1^*\phi_i\cdot c(d)^a$ is equal to 
$n+a-(nd+n+d-1)=-nd+a-d+1$ which is less than zero for $a<d$, hence 
this integral is zero for these $a$. Therefore, 
$$
U_i(d)=\frac{1}{h}\sum_{a\geq d}\frac{(-1)^a}{h^a}
\int_{X_d}e_1^*\phi_i\cdot c(d)^a=
\frac{1}{h^d}\int_{X_d}\frac{e_1^*\phi_i\cdot(-c(d))^d}
{h+c(d)}
$$
Let us denote
$$
u_i(d):=\int_{X_d}\frac{e_1^*\phi_i\cdot(-c(d))^d}
{h+c(d)}
\eqno{(4.9)}
$$
(for $d=0$ we set $u_i(0):=1$). 
Thus, 
$$
z_i(q)=\sum_{d\geq 0}\ u_i(d)q^d;\ u_i(d)=U_i(d)h^d
\eqno{(4.10)}
$$
We have to prove that $z_i(d)=b_i(d)$. 
Therefore, Theorem 2.3 is equivalent to 

{\bf 4.4. Theorem.} {\it The integrals $u_i(d)$ satisfy the relations
$$
u_i(d,h)=\sum_{j\neq i}\sum_{k=1}^d\frac{c_i^j(k)}
{\lambda_i-\lambda_j+kh}\cdot u_j(d-k,\frac{\lambda_j-\lambda_i}{k})
\eqno{(4.10)}
$$}

This is what we are going to prove in the next section, using 
the Bott's localization theorem.

\bigskip
\centerline{\bf \S 5. Fixed point formula.}
\medskip

{\bf 5.1.} We will compute the integrals $u_i(d)$ (see (4.9)) 
by means of the {\it Bott's fixed point formula}. It says that 
$$
u_i(d)=\sum_P\ u_i(P), 
\eqno{(5.1)}
$$
the summation running over all connected components $P\subset X_d^T$. 
Here $u_i(P)$ denotes the integral
$$
u_i(P):=\int_P\left(\frac{e_1^*\phi_i\cdot(-c(d))^d}
{h+c(d)}\right)\biggl|_P\cdot\frac{1}{e(\CN_{P/X_d})}
\eqno{(5.2)}
$$
Here $\CN$ denotes the normal bundle, $e$ the Euler (top Chern) class. 

What do the connected components look like (cf. [K])? 
Let $l_{ij}\subset X$ denote 
the straight line connecting the points $x_i$ and $x_j$. These are the 
curves in $X$ stable under the action of $T$.  

A point in $X_d^T$ is a stable map 
$$
f:(C;y_1,y_2)\lra X
\eqno{(5.3)}
$$
such that $f(y_i)\in X^T=\{x_0,\ldots,x_n\}$ and each irreducible component 
$C_1\subset C$ is mapped either to one of the points $x_i$ --- in this 
case we call $C_1$ {\bf vertical}, or to one of the lines $l_{ij}$ --- in this 
case we call $C_1$ {\bf horizontal}. The map 
$$
f\big|_{C_1}:\ C_1=\BP^1\lra l_{ij}
\eqno{(5.4)}
$$ 
is a finite covering ramified at points $x_i$ and $x_j$. The sum of the 
degrees of these coverings over all horizontal $C_1$ should be equal to 
$d$.  
 
The connected component $P$ to which the point (5.3) belongs, is 
specified by the combinatorial data: which irreducible components of $C$ 
are vertical or horizontal, and the degrees of the coverings (5.4) 
for horizontal components. 

{\bf 5.2.} Let us consider the integral $u_i(P)$ (5.2). Let $f$ as 
in (5.3) be a point in $P$.  
We have 
$$
e^*_1\phi_i\bigl|_P=\phi_i(\lambda_j)
\eqno{(5.5)}
$$
if $f(y_1)=x_j$. Therefore, $u_i(P)$ may be nonzero only if $f(y_1)=x_i$. 
We will suppose this is the case from now on. 

Let $C_1\subset C$ be the irreducible component containing $y_1$.  

{\bf 5.3. Claim.} {\it If $u_i(P)$ is nonzero then $C_1$ is horizontal.}

Indeed, suppose that $C_1$ is vertical. Let us call a {\it 
special point} 
a marked point or a point of intersection of two irreducible components.  

The connected component $P$ has the form 
$\CM_{s+1}\times ?$ where $\CM_{s+1}$ is the 
Deligne-Mumford moduli stack of genus $0$ curves with 
$s+1$ marked points, this whole component mapping 
to $x_i$. A generic curve in $\CM_{s+1}$ contains 
the marked point $y_1$, maybe the marked point $y_2$, 
and $s$ or $s-1$ points of intersection with horizontal 
curves, depending on whether it contains $y_2$ or not; 
$s$ special points altogether.   


We have $\dim(\CM_{s+1})=s-2$. Since the total degree of $f$ is $d$, 
the number of horizontal components does not exceed $d$; therefore,
$s-1\leq d$, hence
$$
\dim(\CM_{s+1})<d
\eqno{(5.6)}
$$
Consequently, $c(d)^d\bigl|_P=0$, i.e. $u_i(P)=0$. 
The claim is proven. $\btu$

{\bf 5.4.} Now let us consider the component $P$ containing  


a stable curve
$$
(f: C=C_1\cup C_2\lra X)\in P\subset X_d^T
\eqno{(5.7)}
$$
where $C_1$ is the irreducible component containing the marked point $y_1$, 
which is mapped with multiplicity $k$ onto the line $l_{ij}$,  
and $C_2$ is all the rest. The map
$$
f_1:=f\bigl|_{C_1}:\ C_1=\BP^1\lra l_{ij}
\eqno{(5.8)}
$$
is the $k$-fold covering ramified only over two points $x_i$ and $x_j$, 
where $1\leq k\leq d$.  
The map
$$
f_2:=f\bigl|_{C_2}:\ C_2\lra X
\eqno{(5.9)}
$$
belongs to a connected component $P_2$ of the fixed point space 
$X_{d-k}^T$. 

We want to compute the integral (5.2). Let us compute the terms under 
the integral. As we have already noted, 
$$
e^*_1\phi_i\bigl|_P=\phi_i(\lambda_i)=e_i=\prod_{b\neq i}(\lambda_i-\lambda_b)
\eqno{(5.10)}
$$
We have 
$$
c(d)\bigl|_P=\frac{\lambda_i-\lambda_j}{k}
\eqno{(5.11)}
$$
{\bf 5.5. Normal bundle.} The most labourous job is to compute the Euler 
class of the normal bundle over $P$. We use the {\it Kontsevich formula}, 
(cf. [K], 3.3.1): the class of $\CN_{P/X_d}$ in the Grothendieck group 
of $T$-equivariant bundles over $P$ is equal to 
$$
[\CN_{P/X_d}]=[H^0(C_1;f^*_1\CT_X)_{0\ at\ y}]-
[H^0(C_1;\CT_{C_1})_{0\ at\ y}]+
$$
$$
+[T_{y;C_1}\otimes T_{y;C_2}]+[T_{y_1;C_1}]+[\CN_{P_2/X_{d-k}}]
\eqno{(5.12)}
$$
(we use the notations describing a bundle by the fiber at a point $f$). 
We have 
$$
[H^0(C_1;f^*_1\CT_X)_{0\ at\ y}]=[H^0(C_1;f^*_1\CT_X)]-
[(f^*_1\CT_X)_{y}]
\eqno{(5.13)}
$$
and
$$
[H^0(C_1;\CT_{C_1})_{0\ at\ y}]=[H^0(C_1;\CT_{C_1})]-
[T_{y;C_1}]
\eqno{(5.14)}
$$
Therefore,
$$
[\CN_{P/X_d}]=([H^0(C_1;f^*_1\CT_X)]-[(f^*_1\CT_X)_{y})])+
$$
$$
+(-[H^0(C_1;\CT_{C_1})]+[T_{y_1;C_1}]+[T_{y;C_1}])+
[T_{y;C_1}\otimes T_{y;C_2}]+[\CN_{P_2/X_{d-k}}]
\eqno{(5.15)}
$$
We have
$$
e([T_{y_1;C_1}])=\frac{\lambda_i-\lambda_j}{k};\ 
e([T_{y;C_1}])=\frac{\lambda_j-\lambda_i}{k}
\eqno{(5.16)}
$$
and
$$
e([H^0(C_1;\CT_{C_1})])=\frac{\lambda_i-\lambda_j}{k}\cdot [0]\cdot
\frac{\lambda_j-\lambda_i}{k}
\eqno{(5.17)}
$$
so that the second bracket in (5.15) gives simply $-[0]$. All zeros in this 
game will cancel out in the final expression for $e(\CN_{P/X_d})$ (see (5.22)  
below)!

{\bf 5.6. Lemma.} {\it We have 
$$
e([H^0(C_1;f^*_1\CT_X)])=\prod_{b=0}^n\prod_{m=0}^k
(\frac{k-m}{k}\lambda_i+\frac{m}{k}\lambda_j-\lambda_b)/[0]
\eqno{(5.18)}
$$}
Note that in the product there two factors equal to $[0]$: they 
correspond to the values $(m,b)=(0,i)$ or $(k,j)$. One of these zeros 
is cancelled. 

{\bf Proof.} We have the exact sequence of vector bundles over $X=\BP^n$:
$$
0\lra\CO_X\lra V^*\otimes\CO_X(1)\lra\CT_X\lra 0
\eqno{(5.19)}
$$
We have 
$$
H^0(C_1;f^*_1(V^*\otimes\CO_X(1)))=H^0(C_1;f^*_1\CO_X(1))
\otimes V^*
$$
The Lie algebra $\ft$ acts on $H^0(C_1;f^*_1\CO_X(1))$ by the characters 
$$
\frac{k-m}{k}\lambda_i+\frac{m}{k}\lambda_j;\ \ m=0,\ldots,k,
$$
and on $V^*$ by the characters $-\lambda_b,\ b=0,\ldots,n$. This 
implies the lemma. $\btu$

The remaining zero in (5.18) will cancel out with the zero from (5.17). 

{\bf 5.7.} We have
$$
e([(f^*_1\CT_X)_{y}])=e([\CT_{x_j;X}])=e_j=\prod_{b\neq j}(\lambda_j
-\lambda_b)
\eqno{(5.20)}
$$ 
Finally, we have 
$$
e([T_{y;C_1}]\otimes [T_{y;C_2}])=\frac{\lambda_j-\lambda_i}{k}+
c(d-k)
\eqno{(5.21)}
$$
Combining everything together, we have proven 

{\bf 5.8. Lemma.} {\it We have 
$$
e(\CN_{P/X_d})=\frac{\prod_{b=0}^n\prod_{m=0;\ (b,m)\neq (i,0),(j,k)}^k
(\frac{k-m}{k}\lambda_i+\frac{m}{k}\lambda_j-\lambda_b)}{e_j}\cdot
$$
$$
\cdot
(\frac{\lambda_j-\lambda_i}{k}+c(d-k))\cdot
e(\CN_{P_2;X_{d-k}})\ \ 
\btu
\eqno{(5.22)}
$$}

{\bf 5.9.} Now we can plug our computations in (5.2). We get
$$
u_i(P)=\frac{1}{k}\cdot\frac{e_i\cdot\left(\frac{\lambda_j-\lambda_i}{k}\right)^d}
{h+\frac{\lambda_i-\lambda_j}{k}}\cdot
$$
$$
\cdot\frac{1}{\prod_{b=0}^n\prod_{m=0;\ (b,m)\neq (i,0),(j,k)}^k
(\frac{k-m}{k}\lambda_i+\frac{m}{k}\lambda_j-\lambda_b)}\cdot
$$
$$
\cdot\int_{P_2}\frac{e_j}{c(d-k)+\frac{\lambda_i-\lambda_j}{k}}\cdot
\frac{1}{e(\CN_{P_2/X_{d-k}})}
\eqno{(5.23)}
$$
The overall factor $1/k$ is due to the group of automorphisms of the 
covering $f_1$, having order $k$. In the big product, the terms 
with $m=0$ give the contribution $e_i$, and the terms with $b=i$ give 
the contribution
$$
\prod_{m=1}^k(-\frac{m}{k}\lambda_i+\frac{m}{k}\lambda_j)=
k!\left(\frac{\lambda_j-\lambda_i}{k}\right)^k
\eqno{(5.24)}
$$
The integral is by definition $U_j(P_2)$ (remember that $P_2$ is the connected 
component of the smaller space $X_{d-k}$!). Thus, we get
$$
u_i(P)=\frac{1}{kh+\lambda_i-\lambda_j}\cdot
$$
$$
\cdot\frac{1}{k!\prod_{b\neq i}\prod_{m=1;\ (b,m)\neq (j,k)}^k
\left(\lambda_i-\lambda_b+\frac{m(\lambda_j-\lambda_i)}{k}\right)}\cdot
\left(\frac{\lambda_j-\lambda_i}{k}\right)^{d-k}\cdot
U_j(P_2)=
$$
$$
=\frac{c_i^j(k)}{kh+\lambda_i-\lambda_j}\cdot
\left(\frac{\lambda_j-\lambda_i}{k}\right)^{d-k}\cdot
U_j(P_2)
\eqno{(5.25)}
$$
This implies formula (4.10) (cf. 4.5). Theorem 4.4, and hence Theorem 2.3 
are proven. $\btu$


\newpage
\bigskip
\centerline{\bf Part II. FLAG SPACES}
\bigskip
\centerline{\bf \S 1. Equivariant cohomology of flag spaces}
\bigskip

{\bf 1.1.} Let $V=\BC^{n+1}$. Let $X$ be the variety of complete flags 
of linear subspaces  
$0\subset V_1\subset\ldots\subset V_n\subset V$, $\dim(V_i)=i$. We set 
$$
D:=\dim(X)=n+(n-1)+\ldots+1=\frac{n(n+1)}{2}
\eqno{(1.1)}
$$  
The torus 
$T=(\BC^*)^{n+1}$ acts on $V$. Namely, if $v_0,\ldots,v_n$ is the standard 
basis of $V$, we put
$$
\alpha\circ\left(\sum z_iv_i\right)=\sum\alpha_iz_iv_i,\ 
\alpha=(\alpha_i)\in T
$$
By functoriality, $T$ acts on $X$. 

We denote by $\ft$ the Lie algebra of $T$, and denote by $\lambda_i\in\ft^*$ 
the character which is the projection on the $i$-th component.   

{\bf 1.2. Fixed points.} Let $F=(V_1\subset \ldots\subset V_n)$ be a flag 
fixed under the action of $T$. Since $TV_1=V_1$, there exists $i_0$ such 
that $V=\BC\cdot e_{i_0}$. Taking the quotient, we get the flag 
$F'=(V_2/V_1\subset\ldots\subset V_n/V_1)$ in the space  
$V'=V/V_1$, fixed under the action of the torus $T'=T/(\BC^*)_{i_0}$.
 By induction on $n$, we conclude that there exists the unique 
permutation $(i_0,\ldots,i_n)$ of the set $\{0,\ldots,n\}$ such that 
$V_k$ is spanned by $e_{i_0},\ldots,e_{i_{k-1}}\ (k=1,\ldots,n)$. 

We shall denote the group of all bijections $w:\ \{0,\ldots,n\}\iso 
\{0,\ldots,n\}$ by $W$, and we shall identify such a bijection with 
the permutation $(w(0),\ldots,w(n))$. 

The previous discussion identifies the set of fixed points $X^T$ with $W$: 
to an element $w\in W$ correponds the flag 
$$
x_{w}=(\BC\cdot e_{w(0)}\subset\ldots\subset\oplus_{i=0}^k\BC
\cdot e_{w(i)}\subset\ldots\subset V)\in X^T
\eqno{(1.2)}
$$

{\bf Tangent spaces.} As a $\ft$-module, the tangent space $T_w:=
T_{X;x_w}$ at the 
point $x_w$, has the set of characters $\{\lambda_{i_p}-\lambda_{i_q},\ 
0\leq p<q\leq n\}$. Hence, its Euler class is given by
$$
e_w:=e(T_w)=\prod_{p<q}\ (\lambda_{i_p}-\lambda_{i_q})
\eqno{(1.3)}
$$   

{\bf Fixed lines.} Given a permutation $w=(i_0,\ldots,i_n)$ and an 
integer $p,\ (1\leq p\leq n)$, let $s_pw$ be the permutation 
with $i_p$ and $i_{p-1}$ transposed, the other entries remaining in place. 

If $x_w=(V_1\subset\ldots\subset V_n)$, let $\ell_{w,s_pw}$ be the projective 
line inside $X$ consisting of all flags of the form
$$
V_1\subset\ldots\subset V_{p-1}\subset V'\subset V_{p+1}\subset\ldots\subset V
$$
where all the subspaces $V_i$ are fixed, and $V'$ is varying. This line 
is fixed under the action of $T$, and contains exactly two fixed points: 
$x_w$ and $x_{s_pw}$.  



However, these are not all fixed lines. In fact, the fixed lines passing 
through one fixed point, correspond to all {\it positive} roots 
(and we have just described the lines corresponding to the {\it simple} 
roots), cf. [BGG].  

Consider the case $w=id$ (identity permutation); so $x_{id}$ 
is the standard flag, with $V_i$ spanned by $\{e_0,\ldots,e_{i-1}\}$. 
Given $p<q$, let $s_{pq}\in W$ be the permutation of $p$ and $q$. 
For a matrix $A=(a_{ij})\in GL(2)$, set 
$$
e_{p;A}=a_{11}e_p+a_{21}e_q;\ e_{q;A}=a_{12}e_p+a_{22}e_q
\eqno{(1.4)}
$$
Let $x_A$ be the flag with the spaces $V_{i;A}$ spanned by $e_{0;A},\ldots, 
e_{i-1;A}$ where $e_{i;A}=e_i$ for $i\neq p,q$. When $A$ runs through 
$GL(2)$, the flags $x_A$ form the projective line $\ell_{e;s_{pq}}=GL(2)/B$,  
stable under the action of $T$, and passing through the fixed 
points $x_e$ and $x_{s_{pq}}$. 

In a similar manner, one defines the $T$-stable line $\ell_{w;s_{pq}w}$ 
passing through $x_w$ and $x_{s_{pq}w}$.

{\bf 1.3. Cohomology.} Let $\CL_i$ be the line bundle over $X$ whose fiber 
over a flag $V_1\subset\ldots\subset V_n$ is equal to $V_{i+1}/V_i$; 
let $u_i\in H^2(X)$ be its first Chern class. The cohomology 
algebra of $X$ is equal to  
$$
H^*(X)=\BC[u_0,\ldots,u_n]/(\sigma_0(u),\ldots,\sigma(u))
\eqno{(1.5)}
$$
Here $\sigma_i(u)$ are the elementary symmetric functions defined by the rule
$$
\prod_{i=0}^n(Z+u_i)=Z^{n+1}+\sigma_0(u)Z^n+\ldots+\sigma_n(u)
\eqno{(1.6)}
$$
The $T$-equivariant cohomology $R=H_T^*(X)$ is an $A=H_T^*(pt)=\BC[\lambda_0,
\ldots,\lambda_n]$-algebra isomorphic to
$$
H^*_T(X)=\BC[u_0,\ldots,u_n;\lambda_0,\ldots,\lambda_n]/
(\sigma_0(u)-\sigma_0(\lambda),\ldots,\sigma_n(u)-\sigma_n(\lambda))
\eqno{(1.7)}
$$
For each $w=(i_0,\ldots,i_n)\in W$, define an element of $H^{2D}_T(X)$ by 
$$
\phi_w=\prod_{p<q}\ (u_{i_p}-\lambda_{i_q})
\eqno{(1.8)}
$$
Let $i_w$ denote the embedding $\{x_w\}\hra X$. 

{\bf 1.4. Lemma.} {\it We have $i^*_w\phi_{w'}=e_w\cdot\delta_{ww'}$.} 

{\bf Proof.} It follows from definitions, that 
$$
i^*_wu_p=\lambda_{w(p)}
\eqno{(1.9)}
$$
On the other hand, 
$$
\prod_{p<q}\ (\lambda_{w(p)}-\lambda_q)=0
\eqno{(1.10)}
$$
for $w\neq e$. This implies the lemma. $\btu$

Set
$$
A'=A[e_w^{-1}]_{w\in W}=A[(\lambda_p-\lambda_q)^{-1}]_{p<q};\ R'=R\otimes_AA'
\eqno{(1.11)}
$$

{\bf 1.5. Corollary.} {\it The elements $\{\phi_w\}_{w\in W}$ form an 
$A'$-basis of the algebra $R'$. The multiplication in $R'$ is recovered 
from the rule
$$
\phi_w\phi_{w'}=e_w\delta_{ww'}
\eqno{(1.12)}
$$}

{\bf Proof.} This is a corollary of the Bott's fixed point theorem, 
cf. I.1.2. $\btu$

\bigskip

\bigskip
\centerline{\bf \S 2. Partition function}
\bigskip

{\bf 2.1.} It is convenient to switch to an arbitrary flag space 
$X=G/B$, $G$ being the simple simply connected algebraic group $G$ 
associated to a finite root system $R$.   
The manifold $X$ is  
acted upon by $T$, the maximal torus of $G$. We identify 
$H_2(X;\BZ)$ with the coroot lattice
$$
H_2(X;\BZ)=Q=\oplus_{i\in I}\ \BZ\cdot\alpha^{\vee}_i,
\eqno{(2.1)}
$$
where $\alpha^{\vee}_i$ are the simple coroots; it contains the submonoid 
of positive coroots
$$
Q^+=\oplus_I\ \BN\cdot\alpha^{\vee}_i\subset Q
\eqno{(2.2)}
$$
The cohomology group $H^2(X;\BZ)$ is the dual weight lattice
$$
H^2(X;\BZ)=P=Char(T)=\oplus_I\ \BZ\cdot\omega_i,
\eqno{(2.3)}
$$
where $\omega_i$ are the fundamental weights; it contains the submonoid 
of dominant weights
$$
P^+=\oplus_I\ \BN\cdot\omega_i\subset P
\eqno{(2.4)}
$$
The root system $R$ lies inside $P$; 
we denote by $R^+\subset R$ the 
subset of positive roots and by $\{\alpha_i\}_{i\in I}$ the set of simple 
roots. 
To each $\lambda\in P=Char(T)$ corresponds the line bundle $\CL_{\lambda}$ 
over $X$, with $c_1(\CL_{\lambda})=\lambda$.  

In the Grothendieck group of $T$-equivariant bundles, the class 
of the tangent bundle of $X$ is equal to the sum
$$
[\CT_X]=\oplus_{\alpha\in R^+}\ [\CL_{\alpha}]
\eqno{(2.5)}
$$
In particular,  
$$
\dim(X)=card(R^+)
\eqno{(2.6)}
$$
Let $W$ be the Weyl group of $R$. The torus $T$ acts on $X$ with 
the finite set of fixed points $\{x_w\}_{w\in W}$, each $x_w$ lying inside 
the corresponding Schubert cell (cf. [BGG]). The fixed lines pass through the  
pairs $x_w,x_{s_{\alpha}w}$ where $s_{\alpha}$ is the reflection 
corresonding to a positive root $\alpha$.

We identify the cohomology ring $A=H^*_T(pt)$ with $\BC[P]=\BC[\alpha_i]_I$. 
The Euler classes of the tangent spaces are given by
$$
e_w=e(T_{x_w;X})=w(e_{id});\ e_{id}=\prod_{\alpha\in R^+}\ \alpha
\eqno{(2.7)}
$$
Set $A'=A[\alpha^{-1}]_{\alpha\in R}$. Bott's Theorem gives the $A'$-base 
$\{\phi_w\}_{w\in W}$ in $H^*_T(X)_{A'}$, such that 
$$
i^*_w\phi_w'=\delta_{ww'}e_w
\eqno{(2.8)}
$$  

{\bf 2.2.} 
To each $\beta\in Q^+$ corresponds the space $X_{\beta}$ of stable 
maps 
$$
f:(C;y_1,y_2)\lra X;\ f_*([C])=\beta;\ g(C)=0
\eqno{(2.9)}
$$
of curves of genus $0$ with two marked points, to $X$. Set 
$|\beta|=\sum\ \nu_i$ if $\beta=\sum\ \nu_i\alpha^{\vee}$.

{\bf 2.3. Lemma.} {\it $\dim(X_{\beta})=2|\beta|+\dim(X)-1$}. 

{\bf Proof.} Let us choose a point (2.9) with $C=\BP^1$. We have 
$$
\dim(X_{\beta})=\dim(T_{f;X_{\beta}})=\dim(H^0(C;f^*\CT_X))-
\dim(H^0(C;\CT_C))+
$$
$$
+\dim(T_{y_1;C})+\dim(T_{y_2;C})
\eqno{(2.10)}
$$
Note that $H^1(C;f^*\CT_X)=0$. Therefore, 
$$
\dim(H^0(C;f^*\CT_X))=\chi(C;f^*\CT_X)=\sum_{\alpha\in R^+}\ 
\chi(C;f^*\CL_{\alpha})=
$$
$$
=\sum_{\alpha\in R^+}\ (\langle\beta,\alpha\rangle+1)=
\langle\beta,2\rho\rangle+card(R^+)=2|\beta|+\dim(X),
\eqno{(2.11)}
$$
since $\langle\alpha_i^{\vee},\rho\rangle=1$ for each $i$. 
Here $\rho$ denotes the half-sum of the positive roots. Obviously, 
$\dim(H^0(C;\CT_C))=3$. Plugging this into (2.10), we get the lemma. 
$\btu$

In (2.11), we have used the equalities
$$
\chi(C;f^*\CL_{\alpha})=\langle\beta,\alpha\rangle+1
\eqno{(2.12)}
$$
This number may be negative if the root system is not simply laced. 

{\bf 2.4.} 
Denote
$$
c(\beta)=c_1(\CT_1)\in H^2_T(X_{\beta})
\eqno{(2.13)}
$$
where $\CT_1$ is the line bundle over $X_{\beta}$ whose fiber over a point 
(2.9) is equal to $T_{y_1;C}$. 

Let $e_1:\ X_{\beta}\lra X$ be the "evaluation at $y_1$" map. Our aim is to 
calculate the formal power series 
$$
Z(q)=\sum_{\beta\in Q^+}\ e_{1*}\bigl(\frac{1}{h+c(\beta)}\bigr)
\cdot q^{\beta}
\eqno{(2.14)}
$$
in variables $q=(q_i)_{i\in I}$, with coefficients in the ring 
$H^*_T(X)(h)$. Here $q^{\beta}=\prod\ q_i^{\nu_i}$, for $\beta=\sum\ \nu_i
\alpha_i^{\vee}$.  

This amounts to calculating $|W|$ power series
$$
Z_w(q)=\langle\phi_w,Z\rangle=\sum_{\beta\in Q^+}\ 
\bigl(\int_{X_{\beta}}\frac{e_1^*\phi_w}{h+c(\beta)}\bigr)\cdot q^{\beta}
\eqno{(2.15)}
$$
{\bf 2.5. Example.} For the root system $A_1$, the series $Z(q)$ is given 
by the expression I, (2.3), with $n=1$. The Weyl group has order $2$; 
there are two series $Z_w(q)$:
$$
Z_{\pm}(q)=\sum_{d=0}^{\infty}\ \frac{q^d}{d!h^d\prod_{m=1}^d(\pm\alpha+mh)}=
\sum_{d=0}^{\infty}\ \frac{q^d}{d^!_{0;h}d^!_{\pm\alpha;h}}
\eqno{(2.16)}
$$
Here we have used the notation
$$
d^!_{\alpha;h}=\prod_{m=1}^d\ (\alpha+mh)
\eqno{(2.17)}
$$
The series $Z_+$ (resp. $Z_-$) corresponds to the trivial (resp. non-trivial) 
element of the Weyl group, cf. I, (2.4) and I.3.6.

\bigskip

\bigskip
\centerline{\bf \S 3. Fixed point computation}
\bigskip

{\bf 3.1.} Let us return to the series (2.15).  
Denote
$$
I_w(\beta)=\int_{X_{\beta}}\frac{e_1^*\phi_w}{h+c(\beta)}
\eqno{(3.1)}
$$ 
Set 
$$
J_w(\beta)=\int_{X_{\beta}}\ \frac{e_1^*\phi_w\cdot(-c(\beta))^{|\beta|}}
{h+c(\beta)}
\eqno{(3.2)}
$$
The next lemma is proved in the same manner as I, (4.10). 

{\bf 3.2. Lemma.} {\it We have 
$$
I_w(\beta)=\frac{1}{h^{|\beta|}}J_w(\beta)\ \btu
$$}

{\bf 3.3.} 
Now, we want to compute the integral $J_w(\beta)$ using the fixed point 
formula. We have
$$
J_w(\beta)=\sum_P\ J_w(P)
\eqno{(3.3)}
$$
where
$$
J_w(P)=\int_P\ \biggl(\frac{e_1^*\phi_w\cdot(-c(\beta))^{|\beta|}}
{h+c(\beta)}\biggr)\biggr|_P\cdot\frac{1}{e(\CN_{P/X_{\beta}})}
\eqno{(3.4)}
$$
Here $P$ denotes a component of the fixed point space $X_{\beta}^T$.
 
Let us compute $J_w(P)$. 
The picture of a connected component is the same as 
in I, \S 5. Thus, let $f:\ (C;y_1,y_2)\lra X$ be a point in 
$P\subset X_{\beta}$. The integral $J_w(P)$ may be non-zero only if 
$f(y_1)=x_w$, so we will assume this.
Let $C=C_1\cup C_2$ where $C_1$ is the connected component containing $y_1$ 
and $C_2$ is all the rest. As in I.5.3, we prove that $C_1$ is horizontal, 
i.e. it covers with some multiplicity $k>0$ a fixed line 
$\ell_{w;s_{\alpha}w}$, for some $\alpha\in R^+$. 


To simplify the notations, assume that $w=id$. We have 
$$
e_1^*\phi_{id}\bigr|_P=e_{id}=e(T_{x_{id};X})=\prod_{\gamma\in R^+}\ \gamma
\eqno{(3.5)}
$$
and 
$$
c(\beta)\bigr|_P=\frac{\alpha}{k}
\eqno{(3.6)}
$$
{\bf 3.4. Lemma. } {\it We have
$$
e(\CN_{P/X_{\beta}})=(-1)^k(k!)^2\cdot\biggl(\frac{\alpha}{k}\biggr)^{2k-1}
\cdot
\prod_{\gamma\in R^+;\ \gamma\neq \alpha}\biggl(
\frac{\prod_{m\geq 1}\ (\gamma-\frac{m}{k}\alpha)}
{\prod_{m\geq 1}\ (s_{\alpha}\gamma-\frac{m}{k}\alpha)}\biggr)\cdot
$$
$$
\cdot\frac{e_{id}}{e_{s_{\alpha}}}\cdot (-\frac{\alpha}{k}+
c(\beta-k\alpha^{\vee}))
$$} 

{\bf Proof.} We have (cf. I.5.5)
$$
[\CN_{P/X_{\beta}}]=[H^0(C_1;f_1^*\CT_X)]-[(f_1^*\CT_X)_y]-
$$
$$
-[H^0(C_1;\CT_{C_1})]+[T_{y_1;C_1}]+[T_{y;C_1}]+
[T_{y;C_1}\otimes T_{y;C_2}]+[\CN_{P_2/X_{\beta-k\alpha^{\vee}}}]
\eqno{(3.7)}
$$
Note that $P_2$ (containing the map $f|_{C_2}$) is the connected component 
of $X_{\beta-k\alpha^{\vee}}$. We have
$$
e([(f_1^*\CT_X)_y])=e([T_{x_{s_{\alpha}};X}])=e_{s_{\alpha}};
\eqno{(3.8)}
$$
next,
$$
e([T_{y_1;C_1}])=\frac{\alpha}{k};\ e([T_{y;C_1}]=-\frac{\alpha}{k},
\eqno{(3.9)}
$$
and
$$
e([H^0(C_1;\CT_{C_1})])=\frac{\alpha}{k}\cdot [0]\cdot -\frac{\alpha}{k}
\eqno{(3.10)}
$$
We have 
$$
e([T_{y;C_1}\otimes T_{y;C_2}])=-\frac{\alpha}{k}+c(\beta-k\alpha^{\vee})
\eqno{(3.11)}
$$

{\bf 3.5. Lemma.} {\it We have 
$$
e([H^0(C_1;f_1^*\CT_X)])=\prod_{\gamma\in R^+}\biggl(
\frac{\prod_{m\geq 0}\ (\gamma-\frac{m}{k}\alpha)}{\prod_{m\geq 1}\  (s_{\alpha}\gamma-\frac{m}{k}\alpha)}\biggr)
\eqno{(3.12)}
$$}

{\bf Proof.} As we have already noted, $[\CT_X]=\sum_{\gamma\in R^+}\ 
[\CL_{\gamma}]$ in the Grothendieck group. Also, $H^1(C_1;f_1^*\CT_X)=0$ 
(convexity of $X$), and our lemma is the corollary of the following one.  

{\bf 3.6. Lemma.} {\it We have 
$$
e([\chi(C_1;f_1^*\CL_{\gamma})])=\frac{\prod_{m\geq 0}\ (\gamma-\frac{m}{k}
\alpha)}{\prod_{m\geq 1}\ (s_{\alpha}\gamma-\frac{m}{k}\alpha)}\ 
\eqno{(3.13)}
$$}

All but finite number 
of terms in the fraction cancel out, and we are left with a finite product 
in the numerator (resp. denominator), if $\langle\gamma,\alpha^{\vee}\rangle$ 
is positive (resp. negative). $\btu$  

Now, Lemma 3.4 follows from Lemma 3.5. Note that the product (3.12) contains 
one factor equal to zero: in the numerator, coresponding to 
$\gamma=\alpha,\ m=k$. This zero cancels out with the zero in (3.10), 
and we are left with the well defined non-zero product. 

This completes the proof of Lemma 3.4. $\btu$ 

Substituting this result into (3.4), we get 

{\bf 3.7. Lemma.} {\it We have
$$
J_w(P;h)=\frac{C_w(\alpha;k)}{kh+w(\alpha)}\cdot 
J_{s_{\alpha}w}(P_2;-\frac{w(\alpha)}{k})
\eqno{(3.14)}
$$
where
$$
C_w(\alpha;k)=w(C_{id}(\alpha;k))
\eqno{(3.15)}
$$
and
$$
C_{id}(\alpha;k)=(-1)^{k(|\alpha|+1)}\alpha^{k|\alpha|-2k+1}\cdot
\frac{k^{k(2-|\alpha|)}}{(k!)^2}\cdot
$$
$$
\cdot\prod_{\gamma\in R^+;\ \gamma\neq \alpha}\biggl(
\frac{\prod_{m\geq 1}\ (s_{\alpha}\gamma-\frac{m}{k}\alpha)}
{\prod_{m\geq 1}\ (\gamma-\frac{m}{k}\alpha)}\biggr)\ \ \btu
\eqno{(3.16)}
$$}

Here we use the notation $|\alpha|=\sum a_i$ for $\alpha=\sum a_i\alpha_i$.

Let us introduce the series
$$
z_w(q)=\sum_{\beta}\ J_w(\beta)q^{\beta}=Z_w(hq)
\eqno{(3.17)}
$$

{\bf 3.8. Theorem.} {\it We have 
$$
z_w(q;h)=1+\sum_{\alpha\in R^+;\ k>0}\ \frac{C_w(\alpha;k)q^{k\alpha^{\vee}}}
{kh+\alpha}\cdot z_{s_{\alpha}w}(q;-\frac{w(\alpha)}{k})
\eqno{(3.18)}
$$
where $C_w(a;k)$ are given by the formulas (3.15), (3.16).}

This is an immediate corollary of Lemma 3.7. $\btu$

Let us look more attentively at the expression (3.16). Recall 

{\bf 3.9. Lemma} ([B], Ch. VI, \S 1, 1.6, Prop. 17, Cor. 2). {\it Let 
$$
w=s_1\cdot\ldots\cdot s_q
\eqno{(3.19)}
$$ 
be a reduced decomposition of an element 
of the Weyl group, where $s_i$ is the reflection correposnding 
to a simple root $\alpha_i$. Then the roots 
$\gamma_i=s_qs_{q-1}\cdot\ldots\cdot s_{i+1}(\alpha_i)\ (i=1,\ldots,q)$ are all 
positive, distinct, and 
$$
R^+\cap w^{-1}(-R^+)=\{\theta_1,\ldots,\theta_q\}\ \ \btu
\eqno{(3.20)} 
$$}

Obviously, if in the product in (3.16) both $\gamma$ and $s_{\alpha}\gamma$ 
are positive, then the two factors corresponding to them cancel out. 
The previous lemma says that we must keep only $l(s_{\alpha})-1$ terms 
in the product (we have already taken care of the term $\gamma=\alpha$).

{\bf 3.10. Corollary.} {\it If the root $\alpha$ is simple then 
$$
C_{id}(\alpha;k)=\frac{k^k}{(k!)^2}\alpha^{-k+1}
\eqno{(3.21)}
$$} with  

In fact, for a simple $\alpha$, the big product disappears at all. 
The expression (3.21) coincides with I, (3.15).


\newpage
\bigskip
\centerline{\bf Part III. COMPUTATIONS FOR $SL(3)$}
\bigskip

\centerline{\bf \S 1. Formula}
\bigskip

{\bf 1.1.} In this part, $X$ will denote the flag manifold $G/B$, 
with $G=SL(3)$.  
According to [G2], the quantum cohomology of $X$ is given by the Fourier 
transform of the following $\CD$-module on the three-torus. 

Let $A=(a_{ij})$ be the $3\times 3$-matrix with 
$a_{ii}=u_{i-1}\ (i=1,2,3)$; $\ a_{i,i-1}=v_i\ (i=2,3)$; 
$\ a_{i,i+1}=-1\ (i=1,2)$; $\ a_{13}=a_{31}=0$. Consider the characteristic polynomial
$$
\det(\lambda+A)=\lambda^3+P_1\lambda^2+P_2\lambda+P_3=
\lambda^3+(u_0+u_1+u_2)\lambda^2+
$$
$$
+(u_0u_1+u_0u_2+u_1u_2+v_1+v_2)\lambda+
u_0u_1u_2+u_0v_2+u_2v_1
\eqno{(1.1)}
$$
Thus, the polynomials $P_i=P_i(u;v)$ are the deformed symmetric functions. 

Consider the three-dimensional torus $T$, with multiplicative coordinates 
$q_0, q_1, q_2$. We define the differential operators $D_i$ on $T$, 
where $D_i$ is obtained from $P_i$ by the substitution $u_j= 
q_j\dpar_{q_j},\ v_j=q_j/q_{j-1}$. We are interested in the solutions 
of the system
$$
D_1\phi=D_2\phi=D_3\phi=0, 
\eqno{(1.2)}
$$
$\phi=\phi(q_0,q_1,q_2)$. 

First of all, since $D_1\phi=0$, the function $\phi$ depends in fact only 
on the quotients $v_1=q_1/q_0,\ v_2=q_2/q_1$. It is useful to write up 
the expressions of the operators $q_i\dpar_{q_i}$ acting on such 
functions, in coordinates $v_1, v_2$: 
$$
q_0\dpar_{q_0}=-v_1\dpar_1;\ 
q_1\dpar_{q_1}=v_1\dpar_1-v_2\dpar_2;\ 
q_2\dpar_{q_2}=v_2\dpar_2
\eqno{(1.3)}
$$
(we set for brevity $\dpar_i=\dpar_{v_i}$).

In these coordinates, the remaining operators look as follows. 
$$
D_2=-(v_1\dpar_1)^2+(v_1\dpar_1)(v_2\dpar_2)-(v_2\dpar_2)^2+v_1+v_2
\eqno{(1.4)}
$$
$$
D_3=-(v_1\dpar_1)^2v_2\dpar_2+v_1\dpar_1(v_2\dpar_2)^2-
v_2(v_1\dpar_1)+v_1(v_2\dpar_2)
\eqno{(1.5)}
$$

{\bf 1.2. Theorem.} {\it There exists a unique, up to a multiplicative 
constant, solution $\phi(v_1,v_2)$ 
of the system 
$$
D_2\phi=D_3\phi=0
\eqno{(1.6)}
$$ 
in the ring of formal power series $\BQ[[v_1,v_2]]$.   

If we normalize $\phi$ by the condition $\phi(0)=1$, it will be given by the formula
$$\
\phi(v_1,v_2)=\sum_{i,j\geq 0}\ \frac{(i+j)!}{(i!)^3(j^!)^3}v_1^iv_2^j
\eqno{(1.7)}
$$}
Here $D_2, D_3$ are given by (1.4), (1.5).

{\bf Proof.} Let us denote the power series (1.7) by 
$\sum_{i,j\geq 0}\ a_{ij}v_1^iv_2^j,\ a_{00}=1$.    
Equation $D_2\phi=0$ is equivalent to the recursion 
$$
(i^2-ij+j^2)a_{ij}=a_{i-1,j}+a_{i,j-1}
\eqno{(1.8)}
$$
(we imply that $a_{ij}=0$ if either $i$ or $j$ is negative). 
Equation $D_3\phi=0$ is equivalent to 
$$
ij(i-j)a_{ij}=-ia_{i,j-1}+ja_{i-1,j}
\eqno{(1.9)}
$$
Both formulas are checked immediately for $a_{ij}$ given by (1.6). 
Already (1.8) defines all $a_{ij}$ uniquely from $a_{00}$. $\btu$

{\bf 1.3. Remarks.} (a) We have $a_{ij}=a_{ji}$. 

(b) It follows from (1.8) that
$$
a_{i0}=\frac{1}{(i!)^2}
\eqno{(1.10)}
$$
The function 
$$
\phi(v,0)=\sum\ \frac{v^i}{(i!)^2}
\eqno{(1.11)}
$$
conicides with the hypergeometric function associated with $G/B$ for 
$G=SL(2)$.

(c) The formulas (1.8) and (1.9) imply the identity
$$
i^3a_{i,j-1}-j^3a_{i-1,j}=0
\eqno{(1.12)}
$$
This, together with (1.10), gives immediately (1.7).    

{\bf 1.4. Another formula.} V. Batyrev communicated to me another formula for 
the solution of (1.6): 
$\tphi=\sum_{i,j}\ b_{ij}v_1^iv_2^j$ where 
$$
b_{ij}=\frac{1}{(i!)^2(j!)^2}\sum_r\ C^r_iC^r_j
\eqno{(1.13)}
$$
where
$$
C^a_b=\frac{b!}{a!(b-a)!}
\eqno{(1.14)}
$$
for integers $a,b$ such that $0\leq a\leq b$, and $0$ if $a<0$ or $a>b$.
 
{\bf Claim.} {\it The function $\tphi$ is equal to $\phi$.} 

Indeed, our claim is equivalent to the identity 
$$
\sum_r\ C^r_iC^r_j=C^i_{i+j}
\eqno{(1.15)}
$$
To prove this, one remarks that to choose $i$ elements from a set 
which is a disjoint union of two sets of cardinalities $j$ and $i$, is 
the same as to choose some $r$ elements from the first set, and 
$i-r$ from the second one.

\bigskip

\newpage
\centerline{\bf \S 2. Equivariant version}

\bigskip

Below, we will deform equations (1.6) and the solution (1.7). 
In terms of quantum cohomology, this deformation corresponds to passing to  
the equivariant cohomology, with respect to a natural action of 
a three-torus on $X$.  

{\bf 2.1.} Introduce differential operators
$$
D_{i;h;\lambda}=P_i(hq_0\dpar_{q_0},hq_1\dpar_{q_1},hq_2\dpar_{q_2};
q_1/q_0,q_2/q_1)-\sigma_i(\lambda_0,\lambda_1,\lambda_2)
\eqno{(2.1)}
$$
Here $P_i(u_0,u_1,u_2;v_1,v_2)$ are the polynomials (1.1); 
$\sigma_i(\lambda)=P_i(\lambda;0)$ are the elementary symmetric functions. 

We are interested in the solutions of the system
$$
D_{1;h;\lambda}\psi=D_{2;h;\lambda}\psi=D_{3;h;\lambda}\phi=0
\eqno{(2.2)}
$$
of the form
$$
\psi(q_0,q_1,q_2)=q_0^{\lambda_0/h}q_1^{\lambda_1/h}q_2^{\lambda_2/h}
\phi(q_0,q_1,q_2);\ \phi(q)=
\sum_{i,j,k}\ b_{ijk}q_0^iq_1^jq_2^k
\eqno{(2.3)}
$$
The equation $D_{1;h;\lambda}\phi=0$ implies $i+j+k=0$, thus, for a solution 
$\psi$, the factor $\phi(q)$ would depend 
only on $v_1=q_1/q_0$ and $v_2=q_2/q_1$. 

To formulate the answer, we introduce the notations $\alpha_i=\lambda_i-
\lambda_{i-1}$; 
$$
p^!_{\alpha}=(h+\alpha)(2h+\alpha)\cdot\ldots\cdot(ph+\alpha)
\eqno{(2.4)}
$$ 
{\bf 2.2. Theorem.} {\it There exists a unique, up to a multiplicative 
constant, solution of the system (2.2) having the form (2.3), 
with $\phi\in\BQ[[v_1,v_2]]$, $v_1=q_1/q_0, v_2=q_2/q_1$. 
If we normalize it by the condition 
$\phi(0)=1$, it will have the form
$$
\phi(v_1,v_2)=\sum_{i,j\geq 0}\ \frac{(i+j)^!_{\alpha_1+\alpha_2}}
{i!j!i^!_{\alpha_1}j^!_{\alpha_2}i^!_{\alpha_1+\alpha_1}j^!_{\alpha_1+\alpha_2}}
\frac{v_1^iv_2^j}{h^{i+j}}
\eqno{(2.5)}
$$}

{\bf Proof.} Let us denote by $a_{ij}$ the coefficient at $v_1^iv_2^j$ 
of the unknown $\phi$. 
The equation $D_{2;h;\lambda}\phi=0$ is equivalent to 
$$
(i^2h^2-ijh^2+j^2h^2+ih\alpha_1+jh\alpha_2)a_{ij}=a_{i-1,j}+a_{i,j-1}
\eqno{(2.6)}
$$
(cf. 1.8)). 
The equation $D_{3;h;\lambda}\phi=0$ is equivalent to 
$$
[(ih-\lambda_0)(ih-jh+\lambda_1)(jh+\lambda_2)+\lambda_0\lambda_1\lambda_2]
a_{ij}=
$$
$$
=(-ih+\lambda_0)a_{i,j-1}+(jh+\lambda_2)a_{i-1,j}
\eqno{(2.7)}
$$
(cf. 1.9)). These identities are checked directly. The uniqueness follows 
from (2.6). $\btu$ 

{\bf 2.3. Remarks.} We have $a_{ij}(\alpha_1,\alpha_2)=
a_{ji}(\alpha_2,\alpha_1)$. Equation (2.6) implies
$$
a_{i0}=\frac{1}{h^ii!i^!_{\alpha_1}}
\eqno{(2.7)}
$$
(cf. (1.10)). 
On the other hand, (2.6) and (2.7) together imply
the identity
$$
ih(ih+\alpha_1)(ih+\alpha_1+\alpha_2)a_{i,j-1}-
jh(jh+\alpha_2)(jh+\alpha_1+\alpha_2)a_{i-1,j}=0
\eqno{(2.8)}
$$
(cf. (1.12)). This and (2.7) give the expression (2.5).

\bigskip
\centerline{\bf \S 3. Comparison with the fixed point method}
\bigskip

{\bf 3.1.} Let us see what does Theorem II.3.8 give for $G=SL(3)$, i.e. 
for the root system $A_2$. The simple roots are $\alpha_1, \alpha_2$, 
the positive ones are $\alpha_1, \alpha_2, \alpha:=\alpha_1+\alpha_2$. 

The Weyl group $W=\Sigma_3$ has two Coxeter generators $s_1, s_2$ 
corresponding to the simple roots, and $s_{\alpha}=s_1s_2s_1$. 
For example, $s_1$ takes $\alpha_1$ to $-\alpha_1$ and $\alpha_2$ to 
$\alpha$, etc. 

We have six series $z_w\ (w\in W)$: 
$$
z_w=z_w(q_1,q_2;\alpha_1,\alpha_2;h)=\sum_{i\geq 0, j\geq 0}\ 
a_{ij;w}q^iq^j;
\eqno{(3.1)}
$$
where
$$
a_{ij;w}=w(a_{ij});\ \ a_{ij}=a_{ij}(\alpha_1,\alpha_2;h):=a_{ij;id}
\eqno{(3.2)}
$$
The element $w$ acts only on the arguments $\alpha_1, \alpha_2$. 

It is easy to see that the recursion relation of Theorem II.3.8 takes 
the following form. 

{\bf 3.2. Theorem.} {\it We have 
$$
z_{id}(q;h)=1+\sum_{k>0}\ \frac{q_1^k}{kh+\alpha_1}\cdot\frac{k^k}{(k!)^2}\cdot 
\frac{1}{\alpha_1^{k-1}}\cdot z_{s_1}(q;-\alpha_1/k)+
$$
$$
+\sum_{k>0}\ \frac{q_2^k}{kh+\alpha_2}\cdot\frac{k^k}{(k!)^2}\cdot 
\frac{1}{\alpha_2^{k-1}}\cdot z_{s_2}(q;-\alpha_2/k)
-\sum_{k>0}\ \biggl(\frac{q_1^kq_2^k}{kh+\alpha_1+\alpha_2}
\cdot\frac{k^{2(k-1)}}{(k!)^2}\cdot
$$
$$
\cdot\frac{\alpha_1+\alpha_2}{\alpha_1\alpha_2}\cdot
\frac{1}{\prod_{m=1}^{k-1}\ (m\alpha_1-(k-m)\alpha_2)^2}\cdot 
z_{s_{\alpha}}(q;-(\alpha_1+\alpha_2)/k)\biggr)\ \ \btu
\eqno{(3.3)}
$$}

Now, the main result of this section is  

{\bf 3.3. Theorem.} {\it The recursion relation} (3.3) {\it is satisfied, with 
$$
a_{ij}=\frac{(i+j)^!_{\alpha_1+\alpha_2}}{i!j!i^!_{\alpha_1}j^!_{\alpha_2}
i^!_{\alpha_1+\alpha_2}j^!_{\alpha_1+\alpha_2}}
\eqno{(3.4)}
$$}

{\bf Proof.} Assume that $i\leq j$. Let us write $a_{ij}$ in the form
$$
a_{ij}(\alpha_1,\alpha_2;h)=\frac{[(i+j)h+\alpha_1+\alpha_2]\cdot 
\ldots\cdot[(1+j)h+\alpha_1+\alpha_2]}
{(jh+\alpha_1+\alpha_2)\cdot\ldots\cdot(h+\alpha_1+\alpha_2)}\cdot
$$
$$
\cdot\frac{1}{i!(h+\alpha_1)\cdot\ldots\cdot(ih+\alpha_1)\cdot
j!(h+\alpha_2)\cdot\ldots\cdot(jh+\alpha_2)}
\eqno{(3.5)}
$$
The denominator (as a function of $h$) has the distinct roots: 
$$
h=-\frac{\alpha_1+\alpha_2}{k}\ (k=1,\ldots,i);\ 
h=-\frac{\alpha_1}{k}\ (k=1,\ldots,i);
$$
$$ 
h=-\frac{\alpha_2}{k}\ (k=1,\ldots,j).
\eqno{(3.6)}
$$
Accordingly, we have the simple fraction decomposition
$$
a_{ij}(\alpha_1,\alpha_2;h)=\sum_{k=1}^i\ 
\frac{^1b^k_{ij}(\alpha_1,\alpha_2)}{kh+\alpha_1}+
\sum_{k=1}^j\ \frac{^2b^k_{ij}(\alpha_1,\alpha_2)}{kh+\alpha_2}+
\sum_{k=1}^i\ \frac{^3b^k_{ij}(\alpha_1,\alpha_2)}{kh+\alpha_1+\alpha_2}
\eqno{(3.7)}
$$
The theorem is equivalent to  

{\bf 3.4. Lemma.} {\it We have} 

(a) {\it For $1\leq k\leq i$,}  
$$
^1b_{ij}^k(\alpha_1,\alpha_2)=\frac{k^k}{(k!)^2}\cdot\frac{1}{\alpha_1^{k-1}} 
\cdot a_{i-k,j}(-\alpha_1,\alpha_1+\alpha_2;-\alpha_1/k);
\eqno{(3.8)}
$$
(b) {\it For $1\leq k\leq j$,}
$$
^2b_{ij}^k(\alpha_1,\alpha_2)=\frac{k^k}{(k!)^2}\cdot\frac{1}{\alpha_2^{k-1}} 
\cdot a_{i,j-k}(\alpha_1+\alpha_2,-\alpha_2;-\alpha_2/k);
\eqno{(3.9)}
$$ 
(c) {\it For $1\leq k\leq i$},
$$
^3b_{ij}^k=-\frac{k^{2(k-1)}}{(k!)^2}\cdot\frac{\alpha_1+\alpha_2}{\alpha_1
\alpha_2}\cdot\frac{1}{\prod_{m=1}^{k-1}\ [m\alpha_1-(k-m)\alpha_2]^2}\cdot
$$
$$
\cdot a_{i-k,j-k}(-\alpha_2,-\alpha_1;-(\alpha_1+\alpha_2)/k)
\eqno{(3.10)}
$$

This lemma is established by a direct computation. The theorem is proved. 
$\btu$

{\bf 3.5. Corollary.} {\it The series $Z_{id}(q_1,q_2)=z_{id}(q_1/h,q_2/h)$ coincides with the series $\phi(q_1,q_2)$ from} (2.5). $\btu$

\bigskip
\centerline{\bf References}
\bigskip

[AB] M.F.~Atiyah and R.~Bott, The moment map and equivariant cohomology, 
{\it Topology} {\bf 23}(1984), 1-28.

[BGG] I.N.~Bernstein, I.M.~Gelfand, S.I.~Gelfand, Schubert cells and 
cohomology of the spaces $G/P$, Usp. Mat. Nauk {\bf 28} (3) (1973), 
3-26 (russian) [Russ. Math. Surv. {\bf 28} (3) (1973), 1-26; 
=I.M.~Gelfand, Coll. papers, Vol. II, 570-595].

[B] N.~Bourbaki, Groupes et alg\`ebres de Lie, Chapitres 4, 5 et 6, Hermann, 
Paris, 1968.   

[G1] A.~Givental, Equivariant Gromov-Witten invariants, {\it IMRN} 
{\bf 13}(1996), 613-663, alg-geom/9603021.

[G2] A.~Givental, Stationary phase integrals, quantum Toda lattices, 
flag manifolds and the Mirror conjecture, in: {\it Topics in Singularity 
Theory, V.I.~Arnold's 60-th Anniversary Collection}, A.~Khovanski\u i, A.~Varchenko, 
V.~Vassiliev (eds.), {\it AMS Translations} Ser. 2, {\bf 180} (1997), 
103-117; alg-geom/9612001.  

[K] M.~Kontsevich, Enumeration of rational curves via torus actions, 
in: {\it The moduli space of curves}, R. Dijkgraaf, C. Faber, van der Geer 
(Eds.), Progress in Math. {\bf 129}, Birkh\"auser, 1995, 335-368, 
hep-th/9405035.

\enddocument